%% file: JP_Oliver.tex
\documentclass[journal]{IEEEtran}
\usepackage{graphicx}
\usepackage{amssymb}
\usepackage[cmex10]{amsmath}
\usepackage{pgfplots, filecontents}
\pgfplotsset{compat=1.12}
\usepackage{epstopdf}
\graphicspath{{./fig/}}
\usepackage{cite}
\usepackage{enumerate}
\usepackage{bm}

\newcommand{\ve}{\mathbf}
\newcommand{\m}{\mathbf}

\newcommand{\vea}[1]{\mathbf{\underline{#1}}}
\newcommand{\ma}[1]{\mathbf{\underline{#1}}}

\newcommand{\re}[1]{\Re\left\{{#1}\right\}} 
\newcommand{\im}[1]{\Im\left\{{#1}\right\}} 

\usepackage{amsthm}
\newtheorem{Proposition}{Result}

\ifCLASSINFOpdf
\else
\fi
\begin{document}
%

\title{Classical Widely Linear Estimation of Real Valued Parameter Vectors in Complex Valued Environments}

%
%
%

\author{Oliver~Lang,~\IEEEmembership{Student Member,~IEEE,}
        and~Mario~Huemer,~\IEEEmembership{Senior Member,~IEEE}
\thanks{O. Lang and M. Huemer are with the Institute of Signal Processing, Johannes Kepler University, Linz, 4040, Austria.}
}
\maketitle

\begin{abstract}
This work investigates the task of estimating a real valued parameter vector based on complex valued measurements in a classical set-up. The application of standard estimators in general results in complex valued estimates of the real valued parameter vector. To avoid this systematic error, widely linear classical estimators that produce real valued estimates are investigated. One of these estimators is the widely linear least squares (WLLS) estimator proposed in this work, which does not utilize any noise statistics. Further, we introduce the best widely linear unbiased estimator (BWLUE) for real valued parameter vectors. The proposed estimators in general outperform their standard counterparts LS estimator and BWLUE, respectively, and they only require half as many complex valued measurements. We compare the novel approaches to standard classical estimators in two application scenarios. One of these applications considers the estimation of a real valued impulse response based on noisy measurements of the system's magnitude and phase response. For this problem, we propose a novel two-step approach based on the introduced widely linear concepts that outperforms standard estimators.
\end{abstract}


%
\IEEEpeerreviewmaketitle

\section{Introduction}
\label{sec:intro}
\input{Introduction}

\section{Preliminaries for Widely Linear Estimators}
\label{sec:AugmentedForm}
\input{AugmentedForm}
\section{BWLUE for Real Valued Parameter Vectors}
\label{sec:Results_improper_Noise}
\input{Results_improper_Noise}

\subsection{Remarks}
\label{sec:Discussion}

\input{Discussion}

\section{Widely Linear Least Squares Estimator for Real Valued Parameter Vectors}
\label{sec:WL-LS}

\input{WL-LS}

\section{Example 1}
\label{sec:Example_1}
\input{Example_1}

\section{Example 2}
\label{sec:Example_2}
\input{Example_2}


\section{Conclusion}
\label{sec:Conclusion}
\input{Conclusion}

\onecolumn

\appendices

\section{Variance and Pseudo Variance of $y_k$}
\label{sec:Variance_of_y_k}
\input{Appendix_A}

\twocolumn

%
%
%
%
%

\ifCLASSOPTIONcaptionsoff
  \newpage
\fi



\bibliography{References}

%

%
%

%

\begin{IEEEbiography}[{\includegraphics[width=1in,height=1.25in,clip,keepaspectratio]{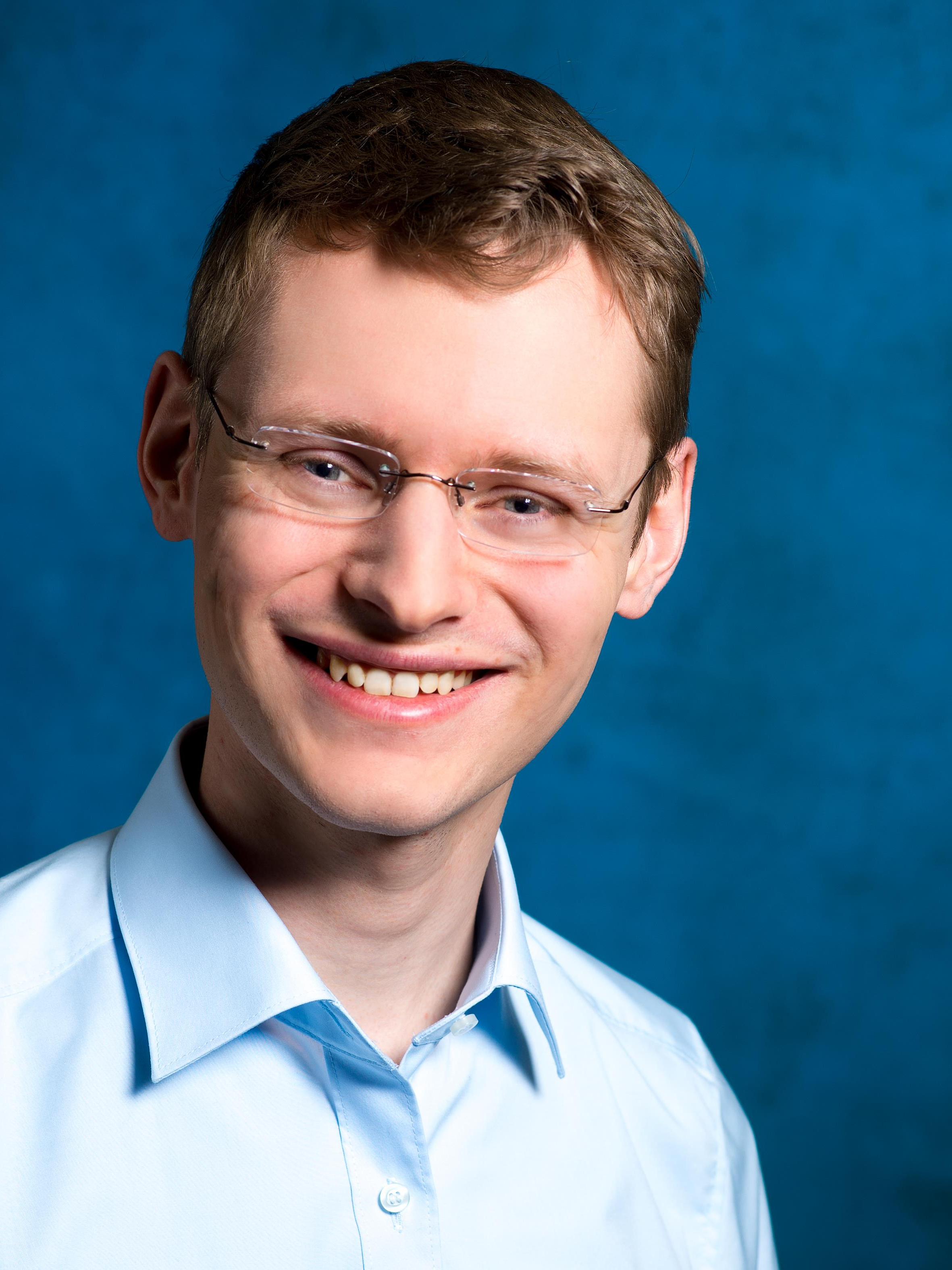}}]{Oliver Lang}
Oliver Lang was born in Sch\"{a}rding, Austria in 1987. He studied Electrical Engineering at the Vienna University of Technology and received his Bachelor degree in 2011. In the next two years Oliver Lang studied Microelectronics at the Vienna University of Technology and finished it with excellent success. The topic of his master thesis was the development and analysis of models for a scanning microwave microscope. Since February 2014, he is a member of the Institute of Signal Processing at the Johannes Kepler University in Linz, Austria.
\end{IEEEbiography}

\begin{IEEEbiography}[{\includegraphics[width=1in,height=1.25in,clip,keepaspectratio]{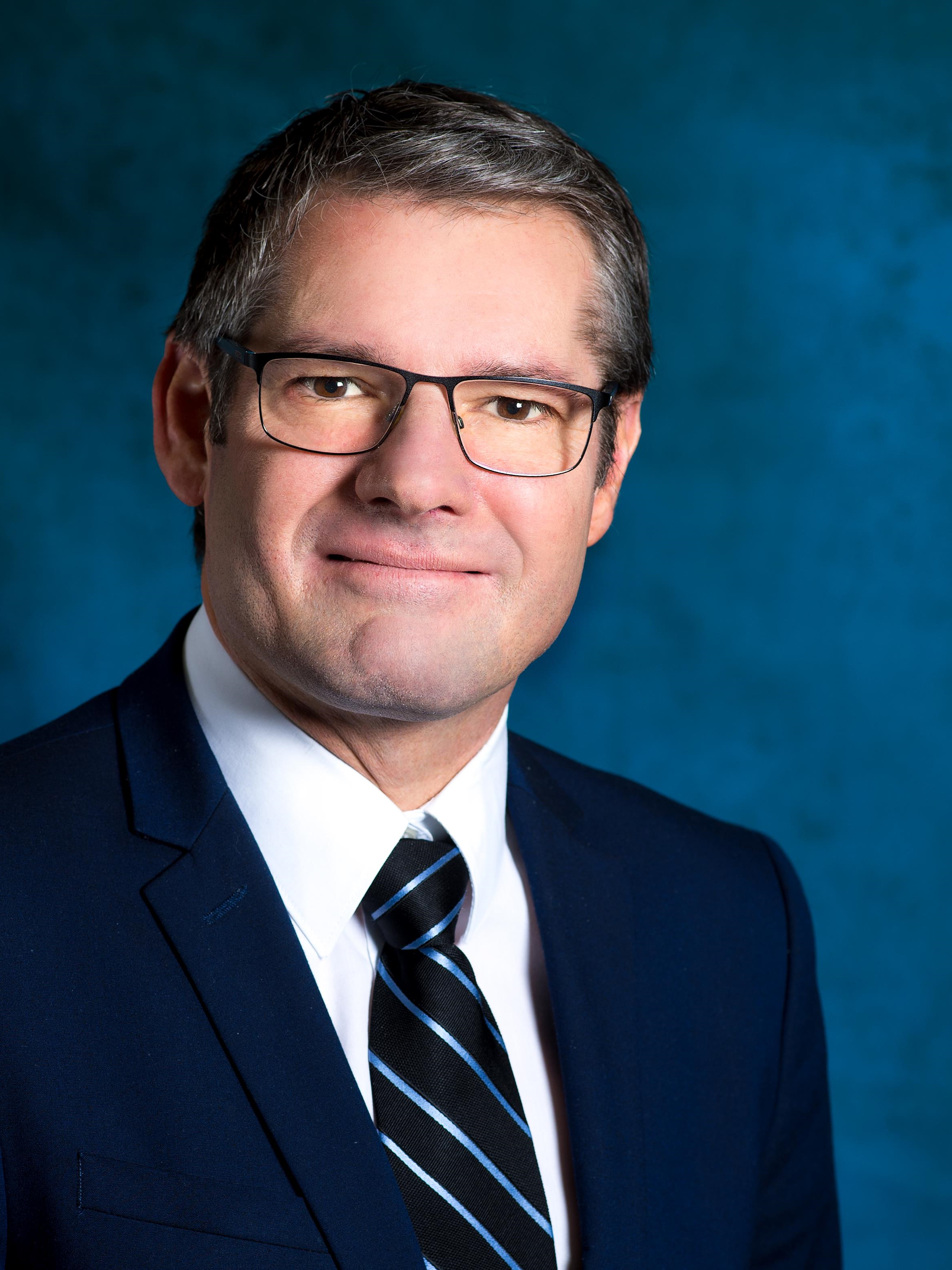}}]{Mario Huemer}
Mario Huemer was born in Wels, Austria, in 1970. He received the Dipl.-Ing. degree in mechatronics and the Dr.techn. (Ph.D.) degree from the Johannes Kepler University (JKU) Linz, Austria, in 1996 and 1999, respectively. From 1997 to 2000, he was a research assistant at the Institute for Communications and Information Engineering at JKU Linz, Austria. From 2000 to 2002, he was with Infineon Technologies Austria, research and development center for wireless products. From 2002-2004 he was a lecturer at the University of Applied Sciences of Upper Austria, and from 2004-2007 he was associate professor for electronics engineering at the University of Erlangen-Nuremberg, Germany. In 2007 Mario Huemer moved to Klagenfurt, Austria, to establish the Chair of Embedded Systems and Signal Processing at Klagenfurt University as a full professor. From 2012 to 2013 he served as dean of the Faculty of Technical Sciences. Since September 2013 he is head of the newly founded Institute of Signal Processing at JKU Linz, Austria.

His research interests are adaptive and statistical signal processing, signal processing architectures and implementations, as well as mixed signal processing with applications in communications, radio frequency and baseband integrated circuits, sensor and
biomedical signal processing. Within these fields he published more than 180 papers. In 2000 Mario Huemer received the German ITG and the Austrian GIT award for dissertations, and in 2010 the Austrian Kardinal Innitzer award in natural sciences. His review work includes national and European research projects as well as international journals. Since 2009 he is member of the editorial board of the "International Journal of Electronics and Communications (AEUE)". 

Mario Huemer is member of the IEEE Signal Processing Society, the IEEE Circuits and Systems Society, and the IEEE Communications Society . He is also member of the German Society of Information Technology  (ITG), and the Austrian Electrotechnical Association (OVE).
\end{IEEEbiography}




\end{document}

%% file: Introduction.tex
%
%
%
%
\IEEEPARstart{I}{n} this work, classical estimation of a real valued parameter vector based on complex valued measurements is investigated. The measurements are assumed to be connected with the parameters via the linear model
\begin{equation}
\ve{y} = \m{H} \ve{x} + \ve{n}, \label{equ:WB_Real001}
\end{equation}
where $\ve{x}\in \mathbb{R}^{N_\ve{x} \times 1}$ is a real valued parameter vector, $\ve{y}\in \mathbb{C}^{N_\ve{y} \times 1}$ is a complex valued measurement vector, $\m{H}\in \mathbb{C}^{N_\ve{y} \times N_\ve{x}}$ is a complex valued measurement matrix, and $\ve{n}\in \mathbb{C}^{N_\ve{y} \times 1}$ is a complex valued zero mean random noise vector. 
A prominent case where such a model appears is the estimation of a real valued impulse response of a linear time-invariant (LTI) system based on complex valued noisy measurements of the system's frequency response. 

In a Bayesian interpretation a real valued parameter vector is improper, and the application of widely linear estimators is obvious \cite{Widely_linear_estimation_with_complex_data, Complex-Valued-Signal-Processing-The-Proper-Way-to-Deal-With-Impropriety, Statistical-Signal-Processing-of-Complex-Valued-Data-The-Theory-of-Improper-and-Noncircular-Signals}. For the definition of propriety we refer to \cite{Complex-Valued-Signal-Processing-The-Proper-Way-to-Deal-With-Impropriety} and Sec.~\ref{sec:AugmentedForm}. A common widely linear Bayesian estimator is the widely linear minimum mean square error (WLMMSE) estimator. The WLMMSE estimator incorporates the fact that $\ve{x}$ is real valued and produces real valued estimates. However, it requires prior knowledge about $\ve{x}$ in form of first and second order statistics. If this kind of prior knowledge is not available, typically classical estimators are employed. In the remainder of this work we only consider classical estimators. 

Note that standard classical estimators such as the linear least squares (LS) estimator \cite{Kay-Est.Theory}, the best linear unbiased estimator (BLUE) \cite{Kay-Est.Theory}, or the best widely linear unbiased estimator (BWLUE) \cite{Statistical-Signal-Processing-of-Complex-Valued-Data-The-Theory-of-Improper-and-Noncircular-Signals} in general do not produce real valued estimates in the setup described above, and hence feature a systematic error. An obvious approach to overcome this issue is to take only the real parts of the estimates for further processing. However, this approach is in general not optimal, as will be further discussed in this work. A special case where this practical approach turns out to be optimal is also discussed. 

Another way to obtain real valued estimates is to rewrite the model in \eqref{equ:WB_Real001} as an equivalent real valued model by using the real composite measurement and noise vector \cite{Complex-Valued-Signal-Processing-The-Proper-Way-to-Deal-With-Impropriety} and apply standard estimators. However, there exist several reasons why the complex notation is preferred in many engineering applications. Among other aspects, the complex representation frequently simplifies the notation, and sometimes gives insights that would not emerge in purely real descriptions \cite{Statistical-Signal-Processing-of-Complex-Valued-Data-The-Theory-of-Improper-and-Noncircular-Signals}. 

Consequently, the first main intend of this work is to propose optimal classical estimators in complex notation incorporating the knowledge that $\ve{x}$ is real valued. We derive an estimator termed BWLUE \emph{for real valued parameter vectors}. It produces real valued estimates, is of widely linear form, and is unbiased in the classical sense, i.e. its estimates $\hat{\ve{x}}$ fulfill $E[\hat{\ve{x}}] = \ve{x}$. The derived estimator considers the general case of improper noise, and allows for a significant complexity reduction if the noise is proper. 

The BWLUE for real valued parameter vectors requires the knowledge of the second order noise statistics, which are not always available. We therefore also derive a widely linear least squares (WLLS) estimator for the model in \eqref{equ:WB_Real001} that does not utilize any noise statistics, and that also produces real valued estimates. This estimator is also extended to the weighted WLLS version. The BWLUE for real valued parameter vectors as well as the proposed WLLS estimator only require half as many complex valued measurements compared to the standard BWLUE and the LS estimator since they only have to estimate half as many real valued parameters. 

The final intend of this work is to apply the derived estimators on practical problems. In a first example we estimate the real valued magnitudes of two complex exponentials based on measurements covered in noise. In a second application we analyze the task of estimating a real valued impulse response of an LTI system based on noisy frequency response measurements. These measurements are often given in form of magnitude and phase response measurements at distinct frequencies. While this is in general a non-linear estimation problem, it can be approximated by a model of the form \eqref{equ:WB_Real001}, however with noise statistics that depend on the unknown parameters to be estimated. To overcome this issue, we discuss several ways to approximate the noise statistics. Furthermore, a novel two-step approach is proposed for this problem, which applies the WLLS estimator proposed in this work in a first step to approximate the noise statistics, and the BWLUE for real valued parameter vectors introduced in this work in a second step.  

The paper is organized as follows: In Sec.~\ref{sec:AugmentedForm}, the basic concepts necessary for investigating widely linear estimators are briefly summarized. In Sec.~\ref{sec:Results_improper_Noise}, the BWLUE for real valued parameter vectors is derived, followed by a number of analytical insights and remarks on this estimator. The WLLS estimator for real valued parameter vectors is proposed in Sec.~\ref{sec:WL-LS}. The subsequent two sections contain the announced applications.  
\\ \\
Notation:
\\ 
Lower-case bold face variables ($\ve{a}$, $\ve{b}$,...) indicate vectors, and upper-case bold face variables ($\m{A}$, $\m{B}$,...) indicate matrices. We further use $\mathbb{R}$ and $\mathbb{C}$ to denote the set of real and complex numbers, respectively, $(\cdot)^T$ to denote transposition, $(\cdot)^H$ to denote conjugate transposition, $(\cdot)^*$ to denote conjugation, $\m{I}^{n\times n}$ to denote the identity matrix of size $n\times n$, and $\m{0}^{m\times n}$ to denote the zero matrix of size $m\times n$. If the dimensions are clear from context we simply write $\m{I}$ and $\m{0}$, respectively. The real and imaginary part of a variable are indicated by $\re{\cdot}$ and $\im{\cdot}$, respectively. 




%% file: AugmentedForm.tex
In this section, we recapitulate the preliminaries required to apply the widely linear estimators discussed in this work. This section is more or less a shortened version of the corresponding parts in \cite{Statistical-Signal-Processing-of-Complex-Valued-Data-The-Theory-of-Improper-and-Noncircular-Signals}.

\subsection{Statistics of Complex-Valued Random Vectors}

We start by constructing the \textit{complex augmented vector} $\underline{\ve{a}}$ of a vector $\ve{a}\in\mathbb{C}^{N_\ve{a}\times 1}$ by stacking $\ve{a}$ on top of its complex conjugate $\ve{a}^{*}$, i.e.
\begin{equation}
	\underline{\ve{a}} = \begin{bmatrix} \ve{a} \\ \ \ve{a}^{*} \end{bmatrix}\in\mathbb{C}^{2N_\ve{a}\times 1}.
\end{equation}
In order to characterize the second-order statistical properties of $\vea{a}$ we consider the augmented covariance matrix 
\begin{align}
	\underline{\m{C}}_{\ve{a}\ve{a}} 
			&= E[(\underline{\ve{a}} - E[\underline{\ve{a}}])(\underline{\ve{a}} - E[\underline{\ve{a}}])^H] \\
			&= \begin{bmatrix} \m{C}_{\ve{a}\ve{a}} & \tilde{\m{C}}_{\ve{a}\ve{a}} \\ 
				 \tilde{\m{C}}_{\ve{a}\ve{a}}^* & \m{C}_{\ve{a}\ve{a}}^* \end{bmatrix} = \underline{\m{C}}_{\ve{a}\ve{a}}^H \in\mathbb{C}^{2N_\ve{a}\times 2N_\ve{a}}, 
				 \label{equ:CWCU_Journal006}
\end{align}
with $\m{C}_{\ve{a}\ve{a}} = E_{\ve{a}}[(\ve{a} - E_{\ve{a}}[\ve{a}])(\ve{a} - E_{\ve{a}}[\ve{a}])^H]$ as the (Hermitian and positive semi-definite) covariance matrix and $\tilde{\m{C}}_{\ve{a}\ve{a}} = E_{\ve{a}}[(\ve{a} - E_{\ve{a}}[\ve{a}])(\ve{a} - E_{\ve{a}}[\ve{a}])^T]$ as the complementary covariance matrix. For $\m{C}_{\ve{a}\ve{a}}$ and $\tilde{\m{C}}_{\ve{a}\ve{a}}$ we have $\m{C}_{\ve{a}\ve{a}} =\m{C}_{\ve{a}\ve{a}}^H$ and $\tilde{\m{C}}_{\ve{a}\ve{a}} = \tilde{\m{C}}_{\ve{a}\ve{a}}^T$, respectively.

In the literature, $\tilde{\m{C}}_{\ve{a}\ve{a}}$ is also referred to as pseudo-covariance matrix or conjugate covariance matrix. If $\tilde{\m{C}}_{\ve{a}\ve{a}} = \m{0}$, then the vector $\ve{a}$ is called \textit{proper}, otherwise 
\textit{improper} \cite{On_the_Circularity_of_a_Complex_Random_Variable, Essential_Statistics_and_Tools_for_Complex_Random_Variables, Second-order_analysis_of_improper_complex_random_vectors_and_processes, A_Complex_Generalized_Gaussian_Distribution_Characterization_Generation_and_Estimation}. 
If $\ve{a}$ is real valued, then $\m{C}_{\ve{a}\ve{a}} = \tilde{\m{C}}_{\ve{a}\ve{a}}$. 

\subsection{Linear and Widely Linear Estimators}

We consider the linear model $\ve{y} = \m{H} \ve{x} + \ve{n}$, where $\ve{x}\in \mathbb{C}^{N_\ve{x} \times 1}$ is an unknown but deterministic parameter vector, $\ve{y}\in \mathbb{C}^{N_\ve{y} \times 1}$ is the measurement vector, $\m{H}\in \mathbb{C}^{N_\ve{y} \times N_\ve{x}}$ is the measurement matrix with full rank and $N_\ve{x} < N_\ve{y}$, and $\ve{n}\in \mathbb{C}^{N_\ve{y} \times 1}$ is a zero mean random noise vector with covariance matrix $\m{C}_{\ve{n}\ve{n}}$ and complementary covariance matrix $\tilde{\m{C}}_{\ve{n}\ve{n}}$. This model can be written in augmented form as
\begin{align}
\vea{y} = \ma{H} \vea{x} + \vea{n},
\end{align}
where
\begin{equation}
\ma{H} = \begin{bmatrix}
\m{H} & \m{0} \\
\m{0} & \m{H}^*
\end{bmatrix}
\end{equation}
and 
\begin{align}
	\underline{\m{C}}_{\ve{n}\ve{n}} 	&=  E[(\underline{\ve{n}} - E[\underline{\ve{n}}])(\underline{\ve{n}} - E[\underline{\ve{n}}])^H] =\begin{bmatrix} \m{C}_{\ve{n}\ve{n}} & \tilde{\m{C}}_{\ve{n}\ve{n}} \\ 
				 \tilde{\m{C}}_{\ve{n}\ve{n}}^* & \m{C}_{\ve{n}\ve{n}}^* \end{bmatrix}.
\end{align}

A widely linear estimator takes on the form 
\begin{equation}
\hat{\ve{x}} = \m{E}\ve{y} + \m{F}\ve{y}^*. \label{equ:WB_Real002}
\end{equation}
In general, widely linear estimators are superior in performance compared to their linear counterparts as soon as the measurements $\ve{y}$ become improper. Applications for widely linear estimators are investigated in \cite{Statistical-Signal-Processing-of-Complex-Valued-Data-The-Theory-of-Improper-and-Noncircular-Signals, Widely-and-semi-widely-linear-processing-of-quaternion-vectors, Optimal-widely-linear-MVDR-beamforming-for-noncircular-signals, Widely_Linear_MVDR_Beamformers_for_the_Reception_of_an_Unknown_Signal_Corrupted_by_Noncircular_Interferences, An_enhanced_widely_linear_CDMA_receiver_with_OQPSK_modulation, On_the_Log-Likelihood_Ratio_Evaluation_of_CWCU_Linear_and_Widely_Linear_MMSE_Data_Estimators, Wise_Conditionally_Unbiased_Widely_Linear_MMSE_Estimation}. Another way to express the estimator in \eqref{equ:WB_Real002} is by its augmented version
\begin{equation}
	\underline{\hat{\ve{x}}} = \begin{bmatrix}  \m{E} & \m{F} \\ \m{F}^{*} & \m{E}^{*} \end{bmatrix} 
														 \begin{bmatrix} \ve{y} \\ \ve{y}^{*} \end{bmatrix}
													 = \underline{\m{G}}\, \underline{\ve{y}},\label{equ:WB_Real004d}
\end{equation}
where $\underline{\m{G}}=\begin{bmatrix}  \m{E} & \m{F} \\ \m{F}^{*} & \m{E}^{*} \end{bmatrix}$. The proposed estimators in this work are compared in performance to the \underline{B}LUE \cite{Kay-Est.Theory} and the \underline{BW}LUE \cite{Statistical-Signal-Processing-of-Complex-Valued-Data-The-Theory-of-Improper-and-Noncircular-Signals} given by
\begin{equation}
\hat{\ve{x}}_\text{B} = \left(\m{H}^H \m{C}_{\ve{n}\ve{n}}^{-1}\,\m{H} \right)^{-1} \m{H}^H \m{C}_{\ve{n}\ve{n}}^{-1}\,\ve{y}, \label{equ:WB_Real004b}
\end{equation}
and
\begin{align}
\hat{\vea{x}}_\text{BW} =& \left(\ma{H}^H \ma{C}_{\ve{n}\ve{n}}^{-1}\,\ma{H} \right)^{-1} \ma{H}^H \ma{C}_{\ve{n}\ve{n}}^{-1}\,\vea{y}, \label{equ:WB_Real004c} 
\end{align}
respectively. Note that the BWLUE reduces to the BLUE for proper noise.

In the following, we discuss some properties of the BWLUE. Considering \eqref{equ:WB_Real004d}, let $\hat{x}_{i}$ denote the $i^\text{th}$ element of $\hat{\ve{x}}$, and let $\ve{e}_{i}^H$,  $\ve{f}_{i}^H$ and $\ve{g}_{i}^H$ denote the $i^\text{th}$ rows of $\m{E}$, $\m{F}$ and $\ma{G}$, respectively. Then, $\hat{x}_{i}$ is given by
\begin{align}
\hat{x}_{i} =& \begin{bmatrix} \ve{e}_{i}^H & \ve{f}_{i}^H  \end{bmatrix} \begin{bmatrix} \ve{y} \\ \ve{y}^*  \end{bmatrix} = \ve{g}_{i}^H \vea{y}, \label{equ:WB_Real007}
\end{align}
with $\ve{g}_{i}^H = \begin{bmatrix} \ve{e}_{i}^H & \ve{f}_{i}^H  \end{bmatrix}$. The BWLUE is defined as the estimator $\hat{x}_{i}$ that minimizes the cost function  \cite{Statistical-Signal-Processing-of-Complex-Valued-Data-The-Theory-of-Improper-and-Noncircular-Signals}
\begin{align}
J =& \text{var}(\hat{x}_{i}) = E\Big[ \left( \hat{x}_{i} - E[\hat{x}_{i}]  \right)\left( \hat{x}_{i} - E[\hat{x}_{i}]  \right)^H   \Big] \\
  =& E\Big[ \left( \ve{g}_{i}^H \vea{y} - E[\ve{g}_{i}^H \vea{y}]  \right)\left( \ve{g}_{i}^H \vea{y} - E[\ve{g}_{i}^H \vea{y}]  \right)^H   \Big] \\
  =& E\Big[ \left(  \ve{g}_{i}^H\vea{n}   \right)\left(  \ve{g}_{i}^H\vea{n} \right)^H   \Big] \\
  =& \ve{g}_{i}^H \ma{C}_{\ve{n}\ve{n}} \ve{g}_{i} \label{equ:WB_Real008}
\end{align}
subject to the unbiased constraint $E[\hat{x}_{i}]=x_i$. From
\begin{align}
E[\hat{x}_{i}] = E[\ve{g}_{i}^H \ma{H} \vea{x} + \ve{g}_{i}^H\vea{n}] = \ve{g}_{i}^H \ma{H} \vea{x} = x_i, \label{equ:WB_Real009c}
\end{align}
it can be seen that the unbiased constraint is fulfilled for every $\vea{x}$ if 
\begin{equation}
\ve{g}_{i}^H \ma{H} = \tilde{\ve{u}}_i^T,\label{equ:WB_Real009}
\end{equation}
where $\tilde{\ve{u}}_i^T$ is a row vector of size $1 \times 2N_\ve{x}$ with a '1' at its $i^\text{th}$ position, and all zeros elsewhere. In summary, the BWLUE for the linear model is the solution of the constrained optimization problem 
\begin{align}
&\ve{g}_{\text{BW},i} = \text{arg} \hspace{2mm} \underset{\ve{g}_i}{\text{min}} \hspace{2mm}\ve{g}_{i}^H \ma{C}_{\ve{n}\ve{n}} \ve{g}_{i}  \hspace{5mm} \text{s.t.} \hspace{3mm} \ve{g}_{i}^H \ma{H} = \tilde{\ve{u}}_i^T, \label{equ:WB_Real010qq}
\end{align}
which can be solved utilizing the Lagrange multiplier method. 

The BLUE and BWLUE in \eqref{equ:WB_Real004b} and \eqref{equ:WB_Real004c} have been derived for complex valued $\ve{x}$. However, in this work we consider real valued $\ve{x}$, while the measurements are considered complex valued. In this case \eqref{equ:WB_Real004c} is no longer the true best widely linear unbiased estimator. The true best widely linear unbiased estimator for real valued $\ve{x}$ will be derived in the next section.

%% file: Results_improper_Noise.tex
\subsection{Derivation of the Estimator} 

In this section, we derive the BWLUE for real valued parameter vectors but complex valued measurements related to each other by the linear model in \eqref{equ:WB_Real001}. In contrast to the ordinary BWLUE in \eqref{equ:WB_Real004c}, the BWLUE for real valued parameter vectors enforces
\begin{align}
\im{\hat{x}_i}& = 0  \label{equ:WB_Real009g} \\
E[\re{\hat{x}_i}] = E[\hat{x}_i]&=x_i. \label{equ:WB_Real009h}
\end{align}
From \eqref{equ:WB_Real009g}, one can easily show that the choice $\ve{e}_i^H = \ve{f}_i^T$ is necessary and sufficient to make $\hat{x}_i$ real valued, independent of the actual realization of $\ve{y}$. Incorporating this into \eqref{equ:WB_Real009h} leads to
\begin{align}
E[ \hat{x}_i ] =& E\left[ \ve{e}_i^H \ve{y} + \ve{e}_i^T \ve{y}^*   \right] \label{equ:WB_Real011} \\
=&\ve{e}_i^H  \m{H} \ve{x} + \ve{e}_i^T \m{H}^* \ve{x} \label{equ:WB_Real013} \\
=& \left(\ve{e}_i^H  \m{H} + \ve{e}_i^T \m{H}^*\right) \ve{x}. \label{equ:WB_Real014} 
\end{align}
Hence, the unbiased constraint $E[ \hat{x}_i ] = \hat{x}_i$ is fulfilled for every $\ve{x}$ if
\begin{equation}
\ve{e}_i^H  \m{H} + \ve{e}_i^T \m{H}^* = \ve{u}_i^T,\label{equ:WB_Real015}
\end{equation}
with $\ve{u}_i^T$ being a row vector of size $1 \times N_\ve{x}$ with a '1' at its $i^\text{th}$ position, and all zeros elsewhere. Together with \eqref{equ:WB_Real008}, we end up at the constrained optimization problem
\begin{align}
\ve{e}_{\mathrm{BW},i} =& \text{arg} \hspace{2mm} \underset{\ve{e}_i}{\text{min}} \hspace{2mm}\left( \begin{bmatrix} \ve{e}_i^H & \ve{e}_i^T  \end{bmatrix} \ma{C}_{\ve{n}\ve{n}} \begin{bmatrix} \ve{e}_i \\ \ve{e}_i^*  \end{bmatrix} \right) \label{equ:WB_Real010ba} \\
=& \text{arg} \hspace{2mm} \underset{\ve{e}_i}{\text{min}} \hspace{2mm} \left( 2 \ve{e}_i^H \m{C}_{\ve{n}\ve{n}} \ve{e}_i + \ve{e}_i^H \widetilde{\m{C}}_{\ve{n}\ve{n}} \ve{e}_i^*   + \ve{e}_i^T  \widetilde{\m{C}}_{\ve{n}\ve{n}}^*  \ve{e}_i \right) \\
& \hspace{5mm} \text{s.t.} \hspace{3mm} \ve{e}_i^H \m{H} + \ve{e}_i^T \m{H}^* = \ve{u}_i^T. \label{equ:WB_Real010}
\end{align}
This can be solved by utilizing the Lagrange multiplier method. The Lagrange cost function follows to
\begin{align}
J' =&  2\ve{e}_i^H \m{C}_{\ve{n}\ve{n}} \ve{e}_i + \ve{e}_i^H \widetilde{\m{C}}_{\ve{n}\ve{n}} \ve{e}_i^* + \ve{e}_i^T  \widetilde{\m{C}}_{\ve{n}\ve{n}}^*  \ve{e}_i  \nonumber \\
& + \bm{\lambda}^T \left( \m{H}^H \ve{e}_i +  \m{H}^T \ve{e}_i^* - \ve{u}_i \right). \label{equ:WB_Real081}
\end{align}
Taking the partial derivative of $J'$ w.r.t $\ve{e}_i^*$ (using Wirtinger's calculus \cite{Wirtinger1927}), results in
\begin{equation}
\frac{\partial J'}{\partial \ve{e}_i^*} =  2 \m{C}_{\ve{n}\ve{n}} \ve{e}_i + 2  \widetilde{\m{C}}_{\ve{n}\ve{n}} \ve{e}_i^*   + \m{H} \bm{\lambda},  \label{equ:WB_Real082}
\end{equation}
where $\bm{\lambda}$ is real valued since the constraint is real valued. Setting the Hermitian of \eqref{equ:WB_Real082} equal to zero and utilizing $\vea{g}_{\mathrm{BW},i}^H = \begin{bmatrix}  \ve{e}_{\mathrm{BW},i}^H & \ve{e}_{\mathrm{BW},i}^T \end{bmatrix}$ yields
\begin{align}
 \ve{e}_{\mathrm{BW},i}^H \m{C}_{\ve{n}\ve{n}} + \ve{e}_{\mathrm{BW},i}^T \widetilde{\m{C}}_{\ve{n}\ve{n}}^* =& - \frac{1}{2}\bm{\lambda}^T \m{H}^H \label{equ:WB_Real083} \\
 \vea{g}_{\mathrm{BW},i}^H \begin{bmatrix}   \m{C}_{\ve{n}\ve{n}} \\  \widetilde{\m{C}}_{\ve{n}\ve{n}}^* \end{bmatrix}    =& - \frac{1}{2}\bm{\lambda}^T \m{H}^H. \label{equ:WB_Real084}
\end{align}
The complex conjugate of \eqref{equ:WB_Real083} can be rewritten in a similar form, producing
\begin{equation}
 \vea{g}_{\mathrm{BW},i}^H \begin{bmatrix}  \widetilde{\m{C}}_{\ve{n}\ve{n}}  \\  \m{C}_{\ve{n}\ve{n}}^* \end{bmatrix}    = - \frac{1}{2}\bm{\lambda}^T \m{H}^T. \label{equ:WB_Real085}
\end{equation}
Combining  \eqref{equ:WB_Real084} and \eqref{equ:WB_Real085} yields
\begin{align}
 \vea{g}_{\mathrm{BW},i}^H  \begin{bmatrix} \m{C}_{\ve{n}\ve{n}} &  \widetilde{\m{C}}_{\ve{n}\ve{n}}  \\ \widetilde{\m{C}}_{\ve{n}\ve{n}}^* & \m{C}_{\ve{n}\ve{n}}^* \end{bmatrix}    =& - \frac{1}{2}\bm{\lambda}^T \underbrace{\begin{bmatrix}  \m{H}^H & \m{H}^T \end{bmatrix}}_{\widetilde{\m{H}}^H} \label{equ:WB_Real086} \\
 \vea{g}_{\mathrm{BW},i}^H \ma{C}_{\ve{n}\ve{n}} =& - \frac{1}{2}\bm{\lambda}^T \widetilde{\m{H}}^H \label{equ:WB_Real087} \\
 \vea{g}_{\mathrm{BW},i}^H =& - \frac{1}{2}\bm{\lambda}^T \widetilde{\m{H}}^H \ma{C}_{\ve{n}\ve{n}}^{-1}, \label{equ:WB_Real088}
\end{align}
where
\begin{align}
\widetilde{\m{H}} = \begin{bmatrix}  \m{H} \\ \m{H}^* \end{bmatrix}. \label{equ:WB_Real100a}
\end{align} 
Inserting \eqref{equ:WB_Real088} into the constraint in \eqref{equ:WB_Real015} produces
\begin{align}
\ve{e}_{\mathrm{BW},i}^H  \m{H} + \ve{e}_{\mathrm{BW},i}^T \m{H}^* =& \ve{u}_i^T  \\
\vea{g}_{\mathrm{BW},i}^H \widetilde{\m{H}} =& \ve{u}_i^T \\
- \frac{1}{2}\bm{\lambda}^T \widetilde{\m{H}}^H \ma{C}_{\ve{n}\ve{n}}^{-1} \widetilde{\m{H}} =& \ve{u}_i^T \\
- \frac{1}{2}\bm{\lambda}^T  =&  \ve{u}_i^T \left( \widetilde{\m{H}}^H \ma{C}_{\ve{n}\ve{n}}^{-1} \widetilde{\m{H}}  \right)^{-1}. \label{equ:WB_Real089}
\end{align}
A reinsertion of \eqref{equ:WB_Real089} into \eqref{equ:WB_Real088} allows to identify $\vea{g}_{\mathrm{BW},i}^H$ as
\begin{align}
 \vea{g}_{\mathrm{BW},i}^H =& \ve{u}_i^T \left( \widetilde{\m{H}}^H \ma{C}_{\ve{n}\ve{n}}^{-1} \widetilde{\m{H}}  \right)^{-1} \widetilde{\m{H}}^H \ma{C}_{\ve{n}\ve{n}}^{-1}. \label{equ:WB_Real090}
\end{align}
The $i^\text{th}$ estimate $\hat{x}_i$ follows to
\begin{align}
 \hat{x}_i =  \ve{e}_{\mathrm{BW},i}^H \ve{y} + \ve{e}_{\mathrm{BW},i}^T \ve{y}^* = \vea{g}_{\mathrm{BW},i}^H \vea{y}. \label{equ:WB_Real100}
\end{align}
Since $\ve{u}_i^T$ is the only term on the right hand side of \eqref{equ:WB_Real090} that depends on the index $i$, the vector estimate $\hat{\ve{x}}$ follows to
\begin{align}
\hat{\ve{x}}_{\mathrm{BW}} =& \left( \widetilde{\m{H}}^H \ma{C}_{\ve{n}\ve{n}}^{-1} \widetilde{\m{H}}  \right)^{-1} \widetilde{\m{H}}^H \ma{C}_{\ve{n}\ve{n}}^{-1} \vea{y} \\
 =& \m{G}_{\mathrm{BW}}\vea{y}, \label{equ:WB_Real101}
\end{align}
where 
\begin{align}
\m{G}_{\mathrm{BW}} = \left( \widetilde{\m{H}}^H \ma{C}_{\ve{n}\ve{n}}^{-1} \widetilde{\m{H}}  \right)^{-1} \widetilde{\m{H}}^H \ma{C}_{\ve{n}\ve{n}}^{-1}. \label{equ:WB_Real102}
\end{align}
Finally, we end up at

\smallskip
\begin{Proposition} \label{prop:Improper_Noise} If $\ve{x}\in\mathbb{R}^{N_\ve{x} \times 1}$ and $\ve{y}\in\mathbb{C}^{N_\ve{y} \times 1}$ are connected via the linear model in \eqref{equ:WB_Real001}, then the BWLUE for real valued parameter vectors is given by $\hat{\ve{x}}_{\mathrm{BW}} = \m{G}_{\mathrm{BW}}\vea{y}$, where the estimator matrix $\m{G}_{\mathrm{BW}}$ is defined in \eqref{equ:WB_Real102} and $\widetilde{\m{H}}$ is defined in \eqref{equ:WB_Real100a}. This estimator is unbiased in the classical sense, i.e. it fulfills $E[\hat{\ve{x}}] = \ve{x}$, and its covariance matrix is
\begin{align}
\m{C}_{\hat{\ve{x}}\hat{\ve{x}},\mathrm{BW}} =& E\left[ \left( \hat{\ve{x}}_{\mathrm{BW}} -E[\hat{\ve{x}}_{\mathrm{BW}}]\right) \left( \hat{\ve{x}}_{\mathrm{BW}} -E[\hat{\ve{x}}_{\mathrm{BW}}]\right) ^H \right] \label{equ:WB_Real115} \\
=& \m{G}_{\mathrm{BW}} \,\ma{C}_{\ve{n}\ve{n}}\, \m{G}_{\mathrm{BW}}^H  \\
=& \left( \widetilde{\m{H}}^H \ma{C}_{\ve{n}\ve{n}}^{-1} \widetilde{\m{H}}  \right)^{-1}. \label{equ:WB_Real116}
\end{align}
\end{Proposition}

\smallskip

The expression for the BWLUE for real valued parameter vectors in Result \ref{prop:Improper_Noise} could have also been derived by minimizing the Bayesian mean square error (BMSE) cost function $E_{\ve{y},\ve{x}}[|\hat{x}_i - x_i|^2], \ \ i=1,\hdots, N_\ve{x}$ subject to the constraint in \eqref{equ:WB_Real015}. This way, the estimator can also be interpreted in a Bayesian sense, where the Bayesian error covariance matrix $\m{C}_{\ve{e}\ve{e}}$ corresponds to $\m{C}_{\hat{\ve{x}}\hat{\ve{x}},\mathrm{BW}}$ in \eqref{equ:WB_Real116} and the minimum BMSEs can be found on the main diagonal of $\m{C}_{\ve{e}\ve{e}}$.

\subsection{Equivalent Real Valued Model}

The complex valued measurements $\ve{y}$ in \eqref{equ:WB_Real001} can also be brought into the form of a real composite vector
\begin{align}
 \ve{y}_\mathbb{R} = \begin{bmatrix}
 \re{\ve{y}} \\ \im{\ve{y}}
\end{bmatrix}.  \label{equ:SPM_040}
\end{align}
$\ve{y}_\mathbb{R}$ is connected with the real valued parameter vector $\ve{x}$ via the real composite linear model
\begin{align}
 \ve{y}_\mathbb{R} &= \underbrace{\begin{bmatrix}
 \re{\m{H}} \\ \im{\m{H}}
\end{bmatrix}}_{\m{H}_\mathbb{R}}\ve{x} + \underbrace{\begin{bmatrix}
 \re{\ve{n}} \\ \im{\ve{n}}
\end{bmatrix}}_{\ve{n}_\mathbb{R}}  \label{equ:SPM_041} \\
&= \m{H}_\mathbb{R}\ve{x} + \ve{n}_\mathbb{R}. \label{equ:WB_Real001a}
\end{align}
For this real valued model, which is equivalent to the complex valued model in \eqref{equ:WB_Real001}, the BLUE minimizing the variances of the elements of $\ve{x}$ subject to the unbiased constraint is given by
\begin{align}
\hat{\ve{x}} =& \left( \m{H}_\mathbb{R}^T \m{C}_{\ve{n}_\mathbb{R}\ve{n}_\mathbb{R}}^{-1} \m{H}_\mathbb{R} \right)^{-1} \m{H}_\mathbb{R}^T \m{C}_{\ve{n}_\mathbb{R}\ve{n}_\mathbb{R}}^{-1} \ve{y}_\mathbb{R},  \label{equ:WB_Real004e}
\end{align}
where $\m{C}_{\ve{n}_\mathbb{R}\ve{n}_\mathbb{R}}$ is connected with $\ma{C}_{\ve{n}\ve{n}}$ according to \cite{Statistical-Signal-Processing-of-Complex-Valued-Data-The-Theory-of-Improper-and-Noncircular-Signals, Complex-Valued-Signal-Processing-The-Proper-Way-to-Deal-With-Impropriety}
\begin{align}
\m{C}_{\ve{n}_\mathbb{R}\ve{n}_\mathbb{R}} = \frac{1}{4}\m{T}^H \ma{C}_{\ve{n}\ve{n}} \m{T}, \label{equ:WB_Real004eqq}
\end{align}
where
\begin{align}
 \m{T} = \begin{bmatrix}
 \m{I}^{N_\ve{y} \times N_\ve{y}} & j  \m{I}^{N_\ve{y} \times N_\ve{y}} \\  \m{I}^{N_\ve{y} \times N_\ve{y}} & -j  \m{I}^{N_\ve{y} \times N_\ve{y}}
 \end{bmatrix}\in \mathbb{C}^{2N_\ve{y} \times 2N_\ve{y}}. \label{equ:WB_Real004q}
\end{align}
By inserting \eqref{equ:WB_Real004eqq} into \eqref{equ:WB_Real004e}, one can easily show that \eqref{equ:WB_Real004e} corresponds to the estimator in complex notation in Result~\ref{prop:Improper_Noise}. However, for the reasons already mentioned in Sec.~\ref{sec:intro} and for several other reasons discussed in \cite{Statistical-Signal-Processing-of-Complex-Valued-Data-The-Theory-of-Improper-and-Noncircular-Signals}, the complex valued representation is often favored.


%% file: Discussion.tex
A discussion and some further properties of the derived estimator are given in this section. We first analyze Result \ref{prop:Improper_Noise} for the case of proper noise. With $\widetilde{\m{C}}_{\ve{n}\ve{n}}= \m{0}^{N_\ve{y} \times N_\ve{y}}$ the estimator in \eqref{equ:WB_Real101}, \eqref{equ:WB_Real102} simplifies to
\begin{align}
\hat{\ve{x}}_{\text{BW}} =&  
\left(\m{H}^H\m{C}_{\ve{n}\ve{n}}^{-1} \m{H} +  \m{H}^T  \left(\m{C}_{\ve{n}\ve{n}}^{-1}\right)^* \m{H}^*  \right)^{-1}   \nonumber \\
&  \cdot \left( \m{H}^H\m{C}_{\ve{n}\ve{n}}^{-1}\ve{y} + \m{H}^T \left(\m{C}_{\ve{n}\ve{n}}^{-1}\right)^*\ve{y}^*  \right)  \\
 = & \left( \re{ \m{H}^H\m{C}_{\ve{n}\ve{n}}^{-1} \m{H}}  \right)^{-1}\re{ \m{H}^H\m{C}_{\ve{n}\ve{n}}^{-1}\ve{y} }. \label{equ:WB_Real036} 
\end{align}
This notation is even simpler as the one for the improper noise case in \eqref{equ:WB_Real101}, \eqref{equ:WB_Real102}, and the evaluation of the estimator becomes significantly less complex.

Assuming the special case, where the term $\m{H}^H\m{C}_{\ve{n}\ve{n}}^{-1} \m{H}$ is real valued we obtain from \eqref{equ:WB_Real036} 
\begin{align}
\hat{\ve{x}}_{\text{BW}} =&    \left(  \m{H}^H\m{C}_{\ve{n}\ve{n}}^{-1} \m{H}  \right)^{-1}\re{ \m{H}^H\m{C}_{\ve{n}\ve{n}}^{-1}\ve{y} } \label{equ:WB_Real038} \\
=& \re{ \left(  \m{H}^H\m{C}_{\ve{n}\ve{n}}^{-1} \m{H}  \right)^{-1}  \m{H}^H\m{C}_{\ve{n}\ve{n}}^{-1}\ve{y}  }. \label{equ:WB_Real039}
\end{align}
In that case, the BWLUE for real valued parameter vectors coincides with the real part of the estimator in \eqref{equ:WB_Real004b}. Furthermore, it also coincides with the real part of the BWLUE in \eqref{equ:WB_Real004c} since the noise is assumed to be proper.

Another interesting statement about the estimator can be made concerning the size of the measurement matrix $\m{H}$. Inspecting \eqref{equ:WB_Real004e} reveals that this estimator is applicable if $\begin{bmatrix} \re{\m{H} } \\ \im{\m{H} } \end{bmatrix} \in \mathbb{R}^{2 N_\ve{y} \times N_\ve{x}}$ has full rank and if $2N_\ve{y} \geq N_\ve{x}$. Therefore, only half as many complex valued measurements are required than there are unknown real valued parameters. This statement clearly also holds for the BWLUE for real valued parameter vectors in Result~\ref{prop:Improper_Noise} since this estimator is equivalent to the one in \eqref{equ:WB_Real004e}.

%% file: WL-LS.tex
The ordinary LS estimator given by
\begin{align}
\hat{\ve{x}} = \left(  \m{H}^H \m{H}  \right)^{-1} \m{H}^H \ve{y}  \label{equ:WB_Real130a} 
\end{align}
can be derived by minimizing the LS cost function \cite{Kay-Est.Theory}
\begin{align}
J(\ve{x}) = \left( \ve{y} - \m{H} \ve{x} \right)^H \left( \ve{y} - \m{H} \ve{x} \right). \label{equ:WB_Real130} 
\end{align}

The expression for the LS estimator in \eqref{equ:WB_Real130a} corresponds to the minimum of the cost function in \eqref{equ:WB_Real130} at least in the following two cases:
\begin{itemize}
\item All terms in \eqref{equ:WB_Real130} are real valued, or
\item all terms in \eqref{equ:WB_Real130} are complex valued.
\end{itemize}
For the case of complex valued measurements but real valued parameters, however, \eqref{equ:WB_Real130a} is no longer optimal in an LS sense. An optimal solution of \eqref{equ:WB_Real130} is derived in this section. It will be shown that this naturally leads to a widely linear estimator which will be termed the WLLS for real valued parameter vectors.  

The expression for the LS estimator in \eqref{equ:WB_Real130a} can also be obtained from the BLUE in \eqref{equ:WB_Real004b} by setting the noise covariance matrix $\m{C}_{\ve{n}\ve{n}}$ equal to the identity matrix (or a scaled version of it). Similarly, it will be shown that the WLLS estimator for real valued parameter vectors also follows from the BWLUE for real valued parameter vectors by setting $\ma{C}_{\ve{n}\ve{n}}$ in Result \ref{prop:Improper_Noise} equal to $\m{I}$. 

The first step for deriving this estimator is to recognize that $J(\ve{x})$ in \eqref{equ:WB_Real130} is real valued even for complex $\ve{y}$ and $\m{H}$. Hence, it can be written in the form 
\begin{align}
&J(\ve{x}) = \nonumber \\
&\frac{1}{2} \left[ \left( \ve{y} - \m{H} \ve{x} \right)^H \left( \ve{y} - \m{H} \ve{x} \right) + \left( \ve{y} - \m{H} \ve{x} \right)^T \left( \ve{y} - \m{H} \ve{x} \right)^* \right]. \label{equ:WB_Real131} 
\end{align}
For real valued $\ve{x}$ but complex $\m{H}$ and $\ve{y}$, the cost function in \eqref{equ:WB_Real131} follows to
\begin{align}
J(\ve{x}) =& \frac{1}{2} \Big[ \ve{y}^H\ve{y} - \ve{y}^H\m{H}\ve{x} - \ve{x}^T \m{H}^H \ve{y} +  \ve{x}^T \m{H}^H \m{H}\ve{x} + 
 \ve{y}^T\ve{y}^* \nonumber \\
 &- \ve{y}^T\m{H}^*\ve{x} - \ve{x}^T \m{H}^T \ve{y}^* +  \ve{x}^T \m{H}^T \m{H}^*\ve{x}  \Big] \label{equ:WB_Real132} \\
 =& \frac{1}{2} \Big[ 2\ve{y}^H\ve{y} - 2\ve{y}^H\m{H}\ve{x} - 2\ve{x}^T \m{H}^H \ve{y} \nonumber \\
 & +  \ve{x}^T \left( \m{H}^H \m{H} +  \m{H}^T \m{H}^*\right)\ve{x}  \Big]. \label{equ:WB_Real133}
\end{align}
Taking the partial derivative of \eqref{equ:WB_Real133} w.r.t. $\ve{x}$ yields
\begin{align}
\frac{\partial J(\ve{x})}{\partial \ve{x}} =& -\m{H}^T \ve{y}^* - \m{H}^H \ve{y} +  \left( \m{H}^H \m{H} +  \m{H}^T \m{H}^*\right)\ve{x} \\
=& - \widetilde{\m{H}}^H \vea{y} + \left( \widetilde{\m{H}}^H \widetilde{\m{H}} \right) \ve{x}. \label{equ:WB_Real134}
\end{align}
Note that no Wirtinger calculus for taking the partial derivative is necessary since $\ve{x}$ is real valued. Setting \eqref{equ:WB_Real134} equal to zero yields
\smallskip
\begin{Proposition} \label{prop:WLLS} If $\ve{x}\in\mathbb{R}^{N_\ve{x} \times 1}$ and $\ve{y}\in\mathbb{C}^{N_\ve{y} \times 1}$ are connected via the linear model in \eqref{equ:WB_Real001}, then the WL\underline{LS} estimator for real valued parameter vectors $\hat{\ve{x}}_{\mathrm{LS}}$ is given by
\begin{equation}
\hat{\ve{x}}_{\mathrm{LS}} = \m{G}_{\mathrm{LS}}\vea{y}, \label{equ:WB_Real061c} 
\end{equation}
where the estimator matrix $\m{G}_{\mathrm{LS}}$ is defined as
\begin{align}
\m{G}_{\mathrm{LS}} =& \left( \widetilde{\m{H}}^H \widetilde{\m{H}} \right)^{-1}  \widetilde{\m{H}}^H. \label{equ:WB_Real061d} 
\end{align}
This result can further be simplified to
\begin{align}
\hat{\ve{x}}_{\mathrm{LS}}  =& \left( \re{\m{H}^H \m{H}} \right)^{-1} \left(\re{ \m{H}^H\ve{y} } \right).  \label{equ:WB_Real061i}
\end{align}
\end{Proposition}
\smallskip
Similar to \eqref{equ:WB_Real038}--\eqref{equ:WB_Real039}, the WLLS estimator reduces to the real part of the LS estimator in \eqref{equ:WB_Real130a} when the term $\m{H}^H \m{H}$ is real valued. Furthermore, the BWLUE from Result~\ref{prop:Improper_Noise} reduces to the WLLS estimator in Result~\ref{prop:WLLS} by setting the augmented noise covariance matrix equal to the identity matrix.

Replacing the LS cost function in \eqref{equ:WB_Real130} by the weighted LS cost function
\begin{align}
J(\ve{x}) = \left( \ve{y} - \m{H} \ve{x} \right)^H \m{W} \left( \ve{y} - \m{H} \ve{x} \right) \label{equ:WB_Real130w} 
\end{align}
allows for deriving the weighted WLLS (WWLLS) estimator for real valued parameter vectors. Assuming that the weighting matrix is a Hermitian matrix makes the derivation a straight forward extension of \eqref{equ:WB_Real131}--\eqref{equ:WB_Real134} and leads to 
\smallskip
\begin{Proposition} \label{prop:WWLLS} If $\ve{x}\in\mathbb{R}^{N_\ve{x} \times 1}$ and $\ve{y}\in\mathbb{C}^{N_\ve{y} \times 1}$ are connected via the linear model in \eqref{equ:WB_Real001}, then the \underline{W}WL\underline{LS} estimator for real parameter vectors $\hat{\ve{x}}_{\mathrm{WLS}}$ is given by
\begin{equation}
\hat{\ve{x}}_{\mathrm{WLS}} = \m{G}_{\mathrm{WLS}}\vea{y} , \label{equ:WB_Real061cw} 
\end{equation}
where the estimator matrix $\m{G}_{\mathrm{WLS}}$ is defined as
\begin{align}
\m{G}_{\mathrm{WLS}} =& \left( \widetilde{\m{H}}^H \ma{W} \widetilde{\m{H}} \right)^{-1}  \widetilde{\m{H}}^H \ma{W}, \label{equ:WB_Real061dw} 
\end{align}
Here, $\ma{W}$ is defined as
\begin{align}
\ma{W} =& \begin{bmatrix} \m{W} & \m{0} \\ \m{0} & \m{W} \end{bmatrix} \label{equ:WB_Real061dwq} 
\end{align}
with $\m{W}$ being a Hermitian weighting matrix. This result can further be simplified to
\begin{align}
\hat{\ve{x}}_{\mathrm{WLS}}  =& \left( \re{\m{H}^H\m{W} \m{H}} \right)^{-1} \left(\re{ \m{H}^H \m{W}\ve{y} } \right).  \label{equ:WB_Real061iw}
\end{align}
\end{Proposition}

\smallskip

Note the similarity between \eqref{equ:WB_Real061dw} and \eqref{equ:WB_Real102}.

%% file: Example_1.tex
In this example, real valued magnitudes of two complex exponentials are estimated based on measurements covered in noise. The measurement at time instance $k$ is written as
\begin{align}
y[k] = x_1 \, \text{exp}\left( j \Omega_1 k \right) + x_2 \, \text{exp}\left( j \Omega_2 k \right) + n[k],
\end{align}
where $k= 1, \cdots, N_\ve{y}$ and where $x_1$ and $x_2$ are the unknown real valued magnitudes. These measurements can be brought into vector/matrix notation
\begin{align}
\ve{y} = \ve{H} \ve{x} + \ve{n}, \label{equ:example1_001}
\end{align}
where $\ve{y}\in \mathbb{C}^{N_\ve{y} \times 1}$ is the measurement vector, $\ve{x}  = \begin{bmatrix} x_1 & x_2 \end{bmatrix}^T$, and
\begin{align}
[\m{H}]_{k,l} = \text{exp}\left( j \Omega_l k \right), \hspace{5mm} l = \{1, 2\}.\label{equ:example1_002}
\end{align}
The noise $\ve{n}$ in \eqref{equ:example1_001} is chosen by \cite{Adaptive_Signal_Processing_Next_Generation_Solutions}
\begin{align}
\ve{n} = \sqrt{1 - \rho^2} \ve{n}_r + j\rho \ve{n}_i, \label{equ:example1_003}
\end{align}
where $\ve{n}_r$ and $\ve{n}_i$ are uncorrelated real valued zero mean Gaussian random vectors of size $N_\ve{y} \times 1$ and with unit variance. With that choice, the noise power remains unaffected while the improperness of the noise can be adjusted by appropriately choosing $\rho$. The noise is proper for $\rho = 1/\sqrt{2}$. In the simulations, we choose $\Omega_1 = 0.1$, $\Omega_2 = 0.2$ and $N_\ve{y}=20$. The following estimators are considered: 
\begin{enumerate}
\item the ordinary LS estimator in \eqref{equ:WB_Real130a},
\item the estimator resulting from taking the real part of the ordinary LS estimator,
\begin{equation}
\hat{\vea{x}} = \re{\left(\m{H}^H \m{H} \right)^{-1} \m{H}^H \,\ve{y}}, \label{equ:WB_Real004cfq}
\end{equation}
\item the WLLS estimator for real valued parameter vectors from Result~\ref{prop:WLLS},
\item the BWLUE in \eqref{equ:WB_Real004c},
\item the estimator resulting from taking the real part of the BWLUE,
\begin{equation}
\hat{\vea{x}} = \re{\left(\ma{H}^H \ma{C}_{\ve{n}\ve{n}}^{-1}\,\ma{H} \right)^{-1} \ma{H}^H \ma{C}_{\ve{n}\ve{n}}^{-1}\,\vea{y}}, \label{equ:WB_Real004cf}
\end{equation}
\item the BWLUE for real valued parameter vectors and improper noise from Result~\ref{prop:Improper_Noise}.
\end{enumerate}

\begin{figure}[tb]
\begin{center}
\begin{tikzpicture}
\begin{semilogyaxis}[compat=newest, 
width=0.9\columnwidth, height = .6\columnwidth, xlabel=$\rho$, 
ylabel style={align=center}, 
ylabel style={text width=3.4cm},
ylabel={average MSE}, 
legend pos=north east, 
legend cell align=left,
legend columns=2, 
        legend style={
            /tikz/column 2/.style={
                column sep=5pt,
            },
        font=\small},
xmin = 0,
xmax = 1,
grid=major,
legend style={
at={(-0.25,1.3)},
anchor=north west}
]

\addplot[line width=1pt, color=red, style=solid] table[x index =0, y index =1] {Average_MSE_x_theor_example_1.dat};
\label{pa3}

\addplot[line width=1pt, color=red, style=dashed] table[x index =0, y index =2] {Average_MSE_x_theor_example_1.dat};
\label{pa5}

\addplot[line width=1pt, color=green, mark repeat=6, mark=triangle*, style=solid, mark options={solid}] table[x index =0, y index =3] {Average_MSE_x_theor_example_1.dat};
\label{pa4}

\addplot[line width=1pt, color=blue] table[x index =0, y index =4] {Average_MSE_x_theor_example_1.dat};
\label{pa1}

\addplot[line width=1pt, color=blue, style=dashed] table[x index =0, y index =5] {Average_MSE_x_theor_example_1.dat};
\label{pa2}

\addplot[line width=1pt, color=green, mark repeat=6, mark=square*] table[x index =0, y index =6] {Average_MSE_x_theor_example_1.dat};
\label{pa6}

\draw [ultra thick, dotted, draw=black] 
        (axis cs: 0.707,0.0000001) -- (axis cs: 0.707,10);
\node at (0.8,0.001) {$1/\sqrt{2}$};

\end{semilogyaxis}
\node [draw,fill=white,anchor=north east] at (rel axis cs: 1.01,2.69) {\shortstack[l]{
	\begin{tabular}{ll}
		\ref{pa3} \small{LS est.} & \ref{pa5} \small{Real part of the LS est.} \\
		\multicolumn{2}{l}{\ref{pa4} \small{WLLS est. for real valued parameter vectors}} \\
  		\ref{pa1} \small{BWLUE}  &
  		\ref{pa2} \small{Real part of the BWLUE} \\
  		\multicolumn{2}{l}{\ref{pa6} \small{BWLUE for real valued parameter vectors }} \\
\end{tabular}}} ;

\end{tikzpicture}
\caption{Average MSEs of the estimated magnitude values for various estimators. The vertical black line marks the value of $\rho = 1/\sqrt{2}$ where the noise is proper. \label{fig:TheorExample1} }
\end{center}
\end{figure}
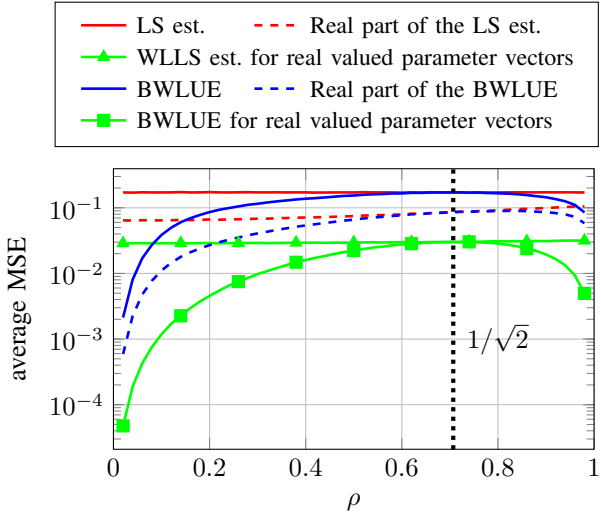

The resulting average mean square errors (MSEs) (averaged over the elements of $\ve{x}$) plotted over $\rho$ are presented in Fig.~\ref{fig:TheorExample1}. The LS estimator performs worst for all values of $\rho$. Its performance can be increased by considering only the real parts of the estimates. Compared to that, a further increase in performance is achieved by the WLLS estimator from Result~\ref{prop:WLLS}. 

The estimators incorporating the improperness of the noise show a performance that strongly depends on $\rho$. One can see that the BWLUE for real valued parameter vectors from Result~\ref{prop:Improper_Noise} outperforms all competing estimators over the whole range of $\rho$. Only for proper noise ($\rho = 1/\sqrt{2}$) the BWLUE from Result~\ref{prop:Improper_Noise} coincides with the WLLS estimator from Result~\ref{prop:WLLS} since $\m{C}_{\ve{n}\ve{n}}$ is a scaled identity matrix and $\tilde{\m{C}}_{\ve{n}\ve{n}}= \m{0}$ in that case.

%% file: Example_2.tex
The second simulation example deals with the old non-linear problem of estimating the sampled impulse response of an analog LTI system based on noisy magnitude and phase response measurements. It will turn out that this example is well suited to test the proposed estimators. Furthermore, it allows to combine the WLLS estimator in Result~\ref{prop:WLLS} with the BWLUE for real valued parameter vectors in Result~\ref{prop:Improper_Noise} to obtain a two-step approach that outperforms the other concepts.

\subsection{Problem Statement} \label{sec:Example_2_Measurement_Model}

The analog real valued impulse response is denoted as $h(t)$. We are interested in estimating the sampled impulse response $h[n] = h(n T_S)$, where $T_S$ is the sampling time. We assume $T_S$ is chosen such that the sampling theorem is practically fulfilled, and we furthermore assume the sampled impulse response to be practically zero after $N_\ve{h}$ samples. Its samples are put together in vector $\ve{h} \in \mathbb{R}^{N_\ve{h} \times 1}$. The measurements are given by $N_\ve{y}$ magnitude and phase response measurements at equidistant frequencies $f_k = k \Delta f$ with $k = 0, \cdots, N_\ve{y}-1$. The true magnitude and phase response values of the analog LTI system at frequency $f_k$ are denoted as $A_k$ and $\varphi_k$, respectively, with $A_k \in \mathbb{R}_0^+$ and $\varphi_k \in [ 0 , 2\pi)$, such that the frequency response $H(f_k)$ is given by  
\begin{align}
H(f_k) = A_k \mathrm{e}^{j \varphi_k}, \hspace{3mm} k = 0, \cdots, N_\ve{y}-1,  
\end{align}
which corresponds to a transformation from polar coordinates to Cartesian coordinates. 
We now define
\begin{align}
H_\mathrm{DC} =& \frac{1}{T_S} H(0) \label{equ:new_007} \\
\ve{H}_\mathrm{AC} =& \frac{1}{T_S} \begin{bmatrix} H(f_1), H(f_1), \hdots, H(f_{N_\ve{y}-1}) \end{bmatrix}^T  \\
\ve{H}_\mathrm{AC,flip} =& \frac{1}{T_S} \begin{bmatrix} H(f_{N_\ve{y}-1}), H(f_{N_\ve{y}-2}), \hdots, H(f_1) \end{bmatrix}^T \label{equ:new_007a}
\end{align}
and
\begin{align}
\ve{H}_\mathrm{ds} = \begin{bmatrix} H_\mathrm{DC} & \ve{H}_\mathrm{AC}^T & \ve{H}_\mathrm{AC,flip}^H \end{bmatrix}^T \in \mathbb{C}^{N_D \times 1}, \label{equ:new_004}
\end{align}
with $N_D = 2 N_\ve{y} -1$. This \underline{d}ouble-\underline{s}ided discrete frequency response is connected with the sampled impulse response according to
\begin{align}
\ve{H}_\mathrm{ds} = \m{F}_\mathrm{ds} \ve{h}, \label{equ:new_001}
\end{align}
where $\m{F}_\mathrm{ds}$ is the matrix given by the first $N_\ve{h}$ columns of the discrete Fourier transform (DFT) matrix of size $N_D \times N_D$. We mainly utilize the \underline{s}ingle-\underline{s}ided frequency response $\ve{H}_\mathrm{ss}$ defined as
\begin{align}
\ve{H}_\mathrm{ss} = \begin{bmatrix} H_\mathrm{DC} & \ve{H}_\mathrm{AC}^T  \end{bmatrix}^T = \m{F}_\mathrm{ss} \ve{h}, \label{equ:new_002}
\end{align}
where $\m{F}_\mathrm{ss}$ is the $N_\ve{y} \times N_\ve{h}$ north west submatrix of the DFT matrix of size $N_D \times N_D$. While the connection between the sampled impulse response $\ve{h}$ and the discrete frequency response in Cartesian coordinates is linear according to \eqref{equ:new_001} or \eqref{equ:new_002}, the relationship between $\ve{h}$ and the magnitude- and phase responses $A_k$ and $\varphi_k$ is non-linear.

\subsection{Measurement Model} \label{sec:Model_Linearization}

We first concentrate on measurements at frequencies $f_k = k \Delta f$ with $k = 1, \cdots, N_\ve{y}-1$, and handle the direct current (DC) measurement later on. The magnitude and phase response measurements at frequency $f_k$ are denoted as $y_k^{(A)}$ and $y_k^{(\varphi)}$, respectively. They are related to $A_k$ and $\varphi_k$ according to
\begin{align}
 y_k^{(A)} &= A_k + n_{A,k} \label{equ:SPM_001} \\
  y_k^{(\varphi)} &= \varphi_k + n_{\varphi,k} \label{equ:SPM_002}
\end{align}
for $k=1,...,N_{\ve{y}}-1$, where $n_{A,k}$ and $n_{\varphi,k}$ denote the corresponding measurement noise variables which we assume to be statistically independent. The probability density functions (PDFs) of $n_{A,k}$ and $n_{\varphi,k}$ clearly depend on the measurement method. We assume $n_{\varphi,k}$ to be zero mean Gaussian with variance $\sigma_{\varphi,k}^2$. Since $A_k$ and $y_k^{(A)}$ have to be positive valued, $n_{A,k}$ cannot be zero mean Gaussian. In our investigations and simulations we consider the following way of generating $y_k^{(A)}$: We sample $n_{A,k}$ from an zero mean Gaussian PDF with variance $\sigma_{A,k}^2$ and add it to $A_k$. If the resulting value of $y_k^{(A)}$ turns out to be negative valued, we set $y_k^{(A)}$ equal to zero. In this way the PDF of $n_{A,k}$ corresponds to a truncated Gaussian with a Delta peak at $n_{A,k} = - A_k$. This PDF of course no longer has zero mean. The mean is denoted as $\mu_k$ in the following. For large $A_k$, the PDF of $n_{A,k}$ is approximately zero mean Gaussian with variance $\sigma_{A,k}^2$, however this is not true for small $A_k$. Transforming the magnitude and phase response measurements to Cartesian coordinates gives
\begin{align}
y_k =& y_k^{(A)} \text{e}^{j y_k^{(\varphi)}} \\
=& (A_k + n_{A,k}) \text{e}^{j (\varphi_k + n_{\varphi,k})}  \label{equ:SPM_003} \\
=& A_k \text{e}^{j \varphi_k} \text{e}^{j n_{\varphi,k}} + n_{A,k} \text{e}^{j \varphi_k} \text{e}^{j n_{\varphi,k}}. \label{equ:SPM_004} 
\end{align}
The random variable $y_k$ can be written as the sum of its mean and a zero mean noise term according to
\begin{align}
y_k =& E\left[ y_k  \right]  + n_{k}.\label{equ:SPM_005}
\end{align}
From \eqref{equ:SPM_004}, the mean $E\left[ y_k  \right]$ follows to
\begin{align}
E\left[ y_k  \right] =& A_k \text{e}^{j \varphi_k} E\left[\text{e}^{j n_{\varphi,k}}  \right] + \mu_k \text{e}^{j \varphi_k} E\left[\text{e}^{j n_{\varphi,k}}  \right] \label{equ:SPM_031}
\end{align}
With $\alpha_k = E\left[\text{e}^{j n_{\varphi,k}}  \right] = E\left[\text{cos}(n_{\varphi,k})\right]= \text{e}^{-\sigma_{\varphi,k}^2 /2} \in [0,1]$ for $n_{\varphi,k} \sim \mathcal{N}(0,\,\sigma_{\varphi,k}^2)$ \cite{Tracking_with_debiased_consistent_converted_measurements_versus_EKF}, and the approximation  $\mu_k \approx 0$ (note that $\mu_k$ depends on the true but unknown magnitude response $A_k$) we have $E\left[ y_k  \right] \approx \alpha_k H(f_k)$ and the measurement model \eqref{equ:SPM_005} for $k = 1, \cdots, N_\ve{y}-1$ simplifies to 
\begin{align}
 y_k \approx \alpha_k H(f_k)  + n_{k}. \label{equ:SPM_031b}
\end{align}
We now turn to the noise term $n_k$ in \eqref{equ:SPM_005}. As shown in Appendix \ref{sec:Variance_of_y_k}, by using the approximation $n_{A,k} \sim \mathcal{N}(0,\,\sigma_{A,k}^2)$ the variance $\sigma_k^2$ and pseudo-variance $\tilde{\sigma}_k^2$ of $n_k$ for $1 \leq k \leq N_\ve{y} -1$ can be approximated by
\begin{align}
\sigma_k^2 =& E[(y_k -E\left[ y_k  \right]) (y_k -E\left[ y_k  \right])^*] \label{equ:SPM_007} \\
\approx & A_k^2\left( 1 - \alpha_k^2 \right) + \sigma_{A,k}^2 \label{equ:SPM_008} 
\end{align}
and
\begin{align}
\tilde{\sigma}_k^2 =& E[(y_k -E\left[ y_k  \right]) (y_k -E\left[ y_k  \right])] \label{equ:SPM_011} \\
\approx & \text{e}^{j 2 \varphi_k} \left( \beta_k A_k^2 + \beta_k \sigma_{A,k}^2 -  A_k^2 \alpha_k^2 \right), \label{equ:SPM_014}
\end{align}
respectively, where  $\beta_k = E\left[\text{e}^{j 2 n_{\varphi,k}}  \right]= E\left[\text{cos}(2 n_{\varphi,k})\right] = \text{e}^{-4 \sigma_{\varphi,k}^2 /2} \in [0,1]$. It is important to note that the noise statistics in \eqref{equ:SPM_008} and \eqref{equ:SPM_014} depend on the true magnitude and phase response values $A_k$ and $\varphi_k$  \cite{Tracking_with_debiased_consistent_converted_measurements_versus_EKF, Unbiased_converted_measurements_for_tracking}. Hence, the true statistics cannot be evaluated without knowing the true magnitude and phase response values. An obvious option is to replace $A_k$ and $\varphi_k$ by $y_k^{(A)}$ and $y_k^{(\varphi)}$ in \eqref{equ:SPM_008} and \eqref{equ:SPM_014}.

\smallskip
We now turn to the measurement at DC, which can be performed by measuring the steady state system response for a unit step at the input. Instead of a magnitude and a phase the measurement at DC is simply given by a real (positive or negative) scalar value denoted by $y_0$. We assume the measurement noise at DC to be zero mean Gaussian with variance $\sigma_0^2 = \sigma_{A,0}^2$ and pseudo-variance $\tilde{\sigma}_0^2 = \sigma_0^2$. 

By defining $y_\mathrm{DC}$, $\ve{y}_\mathrm{AC}$, $\ve{y}_\mathrm{AC,flip}$, $\ve{y}_\mathrm{ds}$ and $\ve{y}_\mathrm{ss}$ according to the rules in \eqref{equ:new_007}--\eqref{equ:new_004} and \eqref{equ:new_002} we finally end up at the compact measurement model 
\begin{align}
	\ve{y}_\mathrm{ss} \approx & T_S \m{D} \m{F}_\mathrm{ss} \ve{h} + \ve{n}, \label{equ:SPM_016} 
\end{align}
where $\m{D}\in \mathbb{R}^{N_{\ve{y}} \times N_{\ve{y}}}$ is a diagonal matrix with $[\m{D}]_{1,1} = 1$ and $[\m{D}]_{k,k} = \alpha_k$ for $k = 1,\cdots, N_\ve{y}-1$.
Assuming the measurements for different $k$ to be statistically independent, the noise covariance matrix $\m{C}_{\ve{n}\ve{n}}$ and the pseudo-noise covariance matrix $\tilde{\m{C}}_{\ve{n}\ve{n}}$ follow to be diagonal matrices that can (according to \eqref{equ:SPM_008} and \eqref{equ:SPM_014}) be approximated by $[\m{C}_{\ve{n}\ve{n}}]_{k,k} = \sigma_k^2$ and $[\tilde{\m{C}}_{\ve{n}\ve{n}}]_{k,k} = \tilde{\sigma}_k^2$ for $k = 0,\cdots, N_\ve{y}-1$. The non-linear connection between the polar measurements and the sampled impulse response has finally been transformed to the model in \eqref{equ:SPM_016} that formally looks like a linear model, but which exhibits noise statistics depending on the true magnitude and phase response values $A_k$ and $\varphi_k$. The noise statistics consequently depend on the unknown vector $\ve{h}$ to be estimated.

\subsection{Estimators} \label{sec:Example_2_Estimators}

In contrast to the first example in Sec.~\ref{sec:Example_1}, we now set the number of measurements $N_{\ve{y}}$ to be smaller than the number of unknown real valued parameters $N_{\ve{h}}$. While this is not an issue for the BWLUE for real valued parameters as discussed in Sec.~\ref{sec:Discussion} as long as $2N_\ve{y} \geq N_\ve{h}$, the ordinary BWLUE fails. Hence, we consider the following estimators:

\begin{enumerate}
\item \textit{IDFT based estimator}: The maybe most intuitive and simple estimator is obtained based on \eqref{equ:new_001} by replacing $\ve{H}_\mathrm{ds}$ with the measurements $\ve{y}_\mathrm{ds}$. An estimate of $\ve{h}$ can be obtained by performing an inverse DFT (IDFT) on $\ve{y}_\mathrm{ds}$ and use the first $N_\ve{h}$ elements of the result as the estimate: 
\begin{align}
\hat{\ve{h}} = \left(\m{F}^{-1}\tilde{\ve{y}}_\mathrm{ds}\right) \odot \ve{w} \label{equ:SPM_018}
\end{align}
Here $\m{F}$ is a DFT matrix of size $N_D \times N_D$ and $\ve{w}\in \mathbb{R}^{N_D \times 1}$ is a windowing vector with ones at the first $N_\ve{h}$ positions and zeros elsewhere. $\odot$ in \eqref{equ:SPM_018} represents the component-wise multiplication. This estimator is in fact a widely linear estimator since it incorporates the measurements and their complex conjugates in a linear way. It always yields a real valued $\ve{h}$, and since it does not incorporate $\m{D}$ it results in biased estimates. 

\item WLLS estimator from Result~\ref{prop:WLLS}: Similar as the IDFT method this estimator does not incorporate the noise statistics. In contrast to the IDFT method the WLLS estimator can also be easily applied when some measurements are missing. In some applications it is for example impossible to measure at DC.

\item BWLUE for real valued parameter vectors from Result~\ref{prop:Improper_Noise}: This estimator is able to incorporate the noise statistics in the form of $\m{C}_{\ve{n}\ve{n}}$ and $\tilde{\m{C}}_{\ve{n}\ve{n}}$. Since in our application the noise statistics depend on the unknown $A_k$ and $\varphi_k$ we insert the measurements $y_k^{(A)}$ and $y_k^{(\varphi)}$ in \eqref{equ:SPM_008} and \eqref{equ:SPM_014} to obtain approximations of the noise statistics. 

\item \textit{Two-step-approach}: Especially when the measurement variances are large, $y_k^{(A)}$ and $y_k^{(\varphi)}$ can deviate heavily from $A_k$ and $\varphi_k$, which might lead to bad approximations of the noise statistics in \eqref{equ:SPM_008} and \eqref{equ:SPM_014}. We therefore suggest the following two-step estimation approach: 1) Perform a WLLS estimation, transform the estimated impulse response into frequency domain using a DFT, and use the resulting frequency response values for approximating the noise statistics in \eqref{equ:SPM_008} and \eqref{equ:SPM_014}. 2) Apply the BWLUE for real valued parameter vectors with the (usually) more precise noise statistics to obtain an improved impulse response estimate.

\end{enumerate}

We have to add one comment to the application of the BWLUE for real valued parameter vectors in this problem: Of course this estimator requires the augmented noise covariance matrix $\ma{C}_{\ve{n}\ve{n}}$ to be invertible. Unfortunately, this is not the case due to $\sigma_0^2 = \tilde{\sigma}_0^2$. However, there exists an easy way to overcome this issue. Consider the real composite model in \eqref{equ:WB_Real001a} with $\m{H}= T_S \m{D} \m{F}_\mathrm{ss}$ (as in \eqref{equ:SPM_016}), and in particular the equation in \eqref{equ:WB_Real001a} corresponding to the first row of $\im{\m{H}}$. First, due to $\sigma_0^2 = \tilde{\sigma}_0^2$ the corresponding diagonal entry of $\m{C}_{\ve{n}_\mathbb{R}\ve{n}_\mathbb{R}}$ is zero (making $\m{C}_{\ve{n}_\mathbb{R}\ve{n}_\mathbb{R}}$ singular). Second, the first row of $\im{\m{H}}$ is a zero row in our problem, such that the according measurement contains no information about $\ve{h}$ at all. As a consequence the corresponding diagonal entry of $\m{C}_{\ve{n}_\mathbb{R}\ve{n}_\mathbb{R}}$ can be set to any arbitrary non-zero value, which makes the noise covariance matrix $\m{C}_{\ve{n}_\mathbb{R}\ve{n}_\mathbb{R}}$ in \eqref{equ:WB_Real004eqq} and consequently also the augmented noise covariance matrix $\ma{C}_{\ve{n}\ve{n}}$ invertible.

\subsection{Simulation Results} 

\begin{figure}[tb]
\begin{center}
\begin{tikzpicture}
\begin{loglogaxis}[compat=newest, 
width=0.9\columnwidth, height = .6\columnwidth, xlabel=$\sigma_{A,k}^2$, 
ylabel style={align=center}, 
ylabel style={text width=3.4cm},
ylabel={average BMSE}, 
legend pos=north east, 
legend cell align=left,
legend columns=2, 
        legend style={
            /tikz/column 2/.style={
                column sep=5pt,
            },
        font=\small},
xmin = 0.0000099,
xmax = 1,
ymax = 0.1,
grid=major,
legend style={
at={(-0.25,1.6)},
anchor=north west}
]

\addplot[line width=1pt, color=black, style=dashed] table[x index =0, y index =5] {Average_BMSE_x_variation_of_C_na.dat};
\label{pc0}

\addplot[line width=1pt, color=green, mark=triangle*, style=solid, mark options={solid}] table[x index =0, y index =2] {Average_BMSE_x_variation_of_C_na.dat};
\label{pc1}

\addplot[line width=1pt, color=green, mark=square*] table[x index =0, y index =1] {Average_BMSE_x_variation_of_C_na.dat};
\label{pc2}

\addplot[line width=1pt, color=red, style=solid] table[x index =0, y index =3] {Average_BMSE_x_variation_of_C_na.dat};
\label{pc3}

\addplot[line width=1pt, color=black, style=solid, style=dotted] table[x index =0, y index =4] {Average_BMSE_x_variation_of_C_na.dat};
\label{pc4}

\end{loglogaxis}
\node [draw,fill=white,anchor=north east] at (rel axis cs: 2.05,3.35) {\shortstack[l]{
	\begin{tabular}{ll}
		\ref{pc0} \small{IDFT method} & \\
		\multicolumn{2}{l}{\ref{pc1} \small{WLLS est. for real valued parameter vectors}} \\
  		\multicolumn{2}{l}{\ref{pc2} \small{BWLUE for real valued parameter vectors}} \\
  		\ref{pc3} \small{Two step approach} & \ref{pc4} \small{Performance bound} 
  	\end{tabular} }};

\end{tikzpicture}
\caption{Average BMSE of the estimated impulse response coefficients for various estimators. The noise variance of the phase response measurements is kept constant at $\sigma_{\varphi,k}^2 = 10^{-1}$ for $1 \leq k \leq N_\ve{y}-1$ while the variances $\sigma_{A,k}^2$ for $0 \leq k \leq N_\ve{y}-1$ are varied between $10^{-5}$ and $1$. \label{fig:TheorExample2} }
\end{center}
\end{figure}
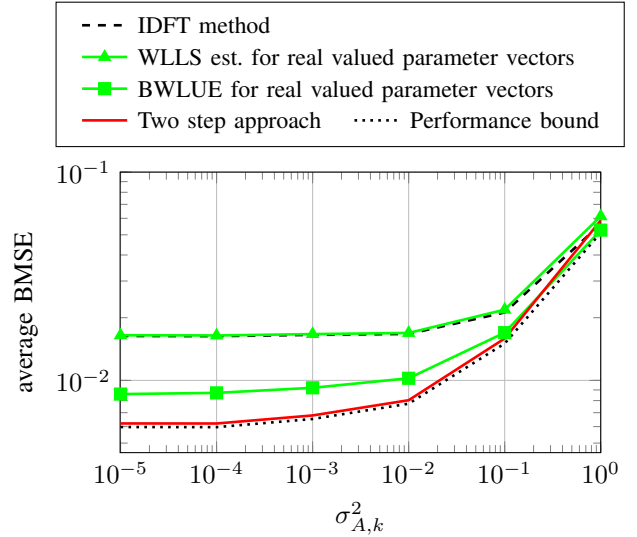

For the simulations, the true impulse responses $\ve{h}$ with length $N_\ve{h} = 12$ are randomly generated by sampling $9$ samples from a normal distribution with zero mean and unit variance, which are then filtered with a finite impulse response (FIR) filter whose coefficients are given by
\begin{align}
\ve{f} = \begin{bmatrix} 0.0881 & 0.4408 & 0.4408 & 0.0881 \end{bmatrix}^T.
\end{align}
This FIR filter corresponds to a low-pass and it takes care that $\ve{h}$ shows low-pass characteristics. Next, the DC- and additional 9 noisy magnitude and phase response measurements are generated as described in Sec. \ref{sec:Model_Linearization}, such that $N_\ve{y} = 10$. In the first experiment the noise variance of the phase response measurements is kept constant at $\sigma_{\varphi,k}^2 = 10^{-1}$ for $1 \leq k \leq N_\ve{y}-1$ while the variances $\sigma_{A,k}^2$ are varied between $10^{-5}$ and $1$ for $0 \leq k \leq N_\ve{y}-1$. Since the true impulse responses are generated randomly, the BMSE is used as a performance measure. The resulting average BMSE curves (averaged over the elements of $\ve{h}$) plotted over $\sigma_{A,k}^2$ are shown in Fig.~\ref{fig:TheorExample2}. In this figure, one can see that the BWLUE for real valued parameter vectors outperforms the WLLS estimator and the IDFT method significantly. By employing the two-step approach, a further increase in performance is achieved. This two-step approach almost reaches our introduced bound, which is simply generated by applying the BWLUE for real valued parameter vectors, but with the true $A_k$ and $\varphi_k$ values inserted in \eqref{equ:SPM_008} and \eqref{equ:SPM_014} to derive the noise statistics.

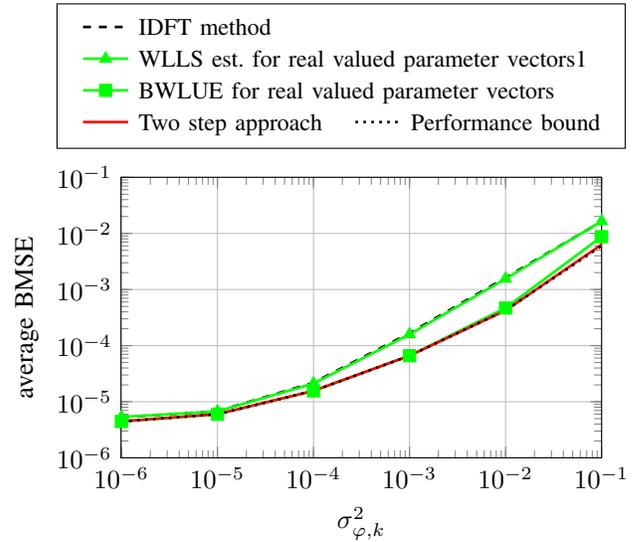
\begin{figure}[tb]
\begin{center}
\begin{tikzpicture}
\begin{loglogaxis}[compat=newest, 
width=0.9\columnwidth, height = .6\columnwidth, xlabel=$\sigma_{\varphi,k}^2$, 
ylabel style={align=center}, 
ylabel style={text width=3.4cm},
ylabel={average BMSE}, 
legend pos=north east, 
legend cell align=left,
legend columns=2, 
        legend style={
            /tikz/column 2/.style={
                column sep=5pt,
            },
        font=\small},
xmin = 0.00000099,
xmax = 0.1,
ymax = 0.1,
ymin = 0.000001,
grid=major,
]

\addplot[line width=1pt, color=black, style=dashed] table[x index =0, y index =5] {Average_BMSE_x_variation_of_C_np.dat};
\label{pb0}

\addplot[line width=1pt, color=green, mark=triangle*, style=solid, mark options={solid}] table[x index =0, y index =2] {Average_BMSE_x_variation_of_C_np.dat};
\label{pb1}

\addplot[line width=1pt, color=green, mark=square*] table[x index =0, y index =1] {Average_BMSE_x_variation_of_C_np.dat};
\label{pb2}

\addplot[line width=1pt, color=red, style=solid] table[x index =0, y index =3] {Average_BMSE_x_variation_of_C_np.dat};
\label{pb3}

\addplot[line width=1pt, color=black, style=solid, style=dotted] table[x index =0, y index =4] {Average_BMSE_x_variation_of_C_np.dat};
\label{pb4}

\end{loglogaxis}
\node [draw,fill=white,anchor=north east] at (rel axis cs: 2.25,2.82) {\shortstack[l]{
	\begin{tabular}{ll}
		\ref{pb0} \small{IDFT method} & \\
		\multicolumn{2}{l}{\ref{pb1} \small{WLLS est. for real valued parameter vectors1}} \\
  		\multicolumn{2}{l}{\ref{pb2} \small{BWLUE for real valued parameter vectors}} \\
  		\ref{pb3} \small{Two step approach} & \ref{pb4} \small{Performance bound} 
  	\end{tabular} }};

\end{tikzpicture}
\caption{Average BMSEs of the estimated impulse response coefficients for various estimators. $\sigma_{A,k}^2 = 10^{-4}$ for $0 \leq k \leq N_\ve{y}-1$ is kept constant while the variances $\sigma_{\varphi,k}^2$ for $1 \leq k \leq N_\ve{y}-1$ are varied between $10^{-6}$ and $10^{-1}$. \label{fig:TheorExample2_b} }
\end{center}
\end{figure}

In the second experiment $\sigma_{A,k}^2 = 10^{-4}$ for $0 \leq k \leq N_\ve{y}-1$ is kept constant while the variances $\sigma_{\varphi,k}^2$ are varied between $10^{-6}$ and $10^{-1}$ for $1 \leq k \leq N_\ve{y}-1$. The resulting average BMSE curves are shown in Fig.~\ref{fig:TheorExample2_b}. This figure shows that the BWLUE for real valued parameter vectors as well as the two-step approach practically reach the bound except for very large values of $\sigma_{\varphi,k}^2$.

%% file: Conclusion.tex
Classical estimation of a real valued parameter vector based on complex valued measurements was investigated in this work. For this task, two widely linear estimators were derived. These estimators are the BWLUE for real valued parameter vectors and the WLLS estimator for real valued parameter vectors. They avoid the systematic error introduced by standard classical estimators by incorporating the fact that the parameter vector is real valued. The proposed estimators were compared with other classical estimators in two application scenarios. For the problem of estimating a real valued impulse response based on noisy frequency response measurements a novel two-step estimation approach was proposed which is adapted from the introduced widely linear estimators in this work.

%% file: Appendix_A.tex
In the following, we make the approximation $n_{A,k} \sim \mathcal{N}(0,\,\sigma_{A,k}^2)$ for $k = 1,\hdots, N_\ve{y}-1$. The variance $\sigma_k^2$ of the $k^\text{th}$ measurement $y_k$ in Cartesian coordinates can be derived as
\begin{align}
\sigma_k^2 =& E\left[\left( y_k - E\left[ y_k  \right] \right)  \left( y_k - E\left[ y_k  \right] \right)^* \right] \\
=& E\left[\left( A_k \text{e}^{j \varphi_k} \text{e}^{j n_{\varphi,k}} + n_{A,k} \text{e}^{j \varphi_k} \text{e}^{j n_{\varphi,k}} - \alpha_k A_k \text{e}^{j \varphi_k} \right)\left( A_k \text{e}^{-j \varphi_k} \text{e}^{-j n_{\varphi,k}} + n_{A,k} \text{e}^{-j \varphi_k} \text{e}^{-j n_{\varphi,k}} - \alpha_k A_k \text{e}^{-j \varphi_k} \right) \right] \\
=& E\left[
A_k^2  + A_k n_{A,k} - \alpha_k A_k^2 \text{e}^{j n_{\varphi,k}} + 
A_k n_{A,k} + n_{A,k}^2 - \alpha_k A_k n_{A,k} \text{e}^{j n_{\varphi,k}} -
\alpha_k A_k^2 \text{e}^{-j n_{\varphi,k}} - \alpha_k A_k n_{A,k} \text{e}^{-j n_{\varphi,k}} +
\alpha_k^2  A_k^2 
\right] \\
=& A_k^2  - \alpha_k^2 A_k^2  
+ \sigma_{A,k}^2 - \alpha_k^2 A_k^2 + \alpha_k^2  A_k^2 \\
=& A_k^2(1 - \alpha_k^2) + \sigma_{A,k}^2 .
\end{align}

Similarly, the pseudo-variance $\tilde{\sigma}_k^2$ of the $k^\text{th}$ measurement $y_k$ in Cartesian coordinates follows as
\begin{align}
\tilde{\sigma}_k^2 =& E\left[\left( y_k - E\left[ y_k  \right] \right)  \left( y_k - E\left[ y_k  \right] \right) \right] \\
=& E\left[\left( A_k \text{e}^{j \varphi_k} \text{e}^{j n_{\varphi,k}} + n_{A,k} \text{e}^{j \varphi_k} \text{e}^{j n_{\varphi,k}} - \alpha_k A_k \text{e}^{j \varphi_k} \right)\left( A_k \text{e}^{j \varphi_k} \text{e}^{j n_{\varphi,k}} + n_{A,k} \text{e}^{j \varphi_k} \text{e}^{j n_{\varphi,k}} - \alpha_k A_k \text{e}^{j \varphi_k} \right) \right] \\
=& E \Big[
 A_k^2 \text{e}^{j2 \varphi_k} \text{e}^{j2 n_{\varphi,k}} + 2 A_k n_{A,k} \text{e}^{j2 \varphi_k} \text{e}^{j2 n_{\varphi,k}}
- 2 \alpha_k A_k^2 \text{e}^{j2 \varphi_k} \text{e}^{j n_{\varphi,k}}  \nonumber \\
& \hspace{5mm} + n_{A,k}^2 \text{e}^{j2\varphi_k} \text{e}^{j2 n_{\varphi,k}}
- 2 \alpha_k A_k n_{A,k} \text{e}^{j2 \varphi_k} \text{e}^{j n_{\varphi,k}} + \alpha_k^2 A_k^2 \text{e}^{j 2\varphi_k}
\Big] \\
=& A_k^2 \text{e}^{j2 \varphi_k} \underbrace{E[\text{e}^{j2 n_{\varphi,k}}]}_{\beta_k} 
- 2 \alpha_k^2 A_k^2 \text{e}^{j2 \varphi_k}  + \sigma_{A,k}^2\text{e}^{j2\varphi_k} E[\text{e}^{j2 n_{\varphi,k}}]
+ \alpha_k^2 A_k^2 \text{e}^{j 2\varphi_k} \\
=& A_k^2 \beta_k \text{e}^{j2 \varphi_k} 
-  \alpha_k^2 A_k^2 \text{e}^{j2 \varphi_k}  + \sigma_{A,k}^2 \beta_k \text{e}^{j2\varphi_k} \\
=& \text{e}^{j2 \varphi_k}  \left( \beta_k A_k^2  - \alpha_k^2 A_k^2   + \sigma_{A,k}^2 \beta_k \right).
\end{align}